\documentclass[12pt,reqno]{amsart}
\usepackage{amsmath}
\usepackage{amssymb}

\def\pg{$p$-group}
\def\agemo{\mho}

\newtheorem{thm}{Theorem}
\newtheorem{lem}{Lemma}

\newtheorem{thmA}{Theorem A (\cite[III.14.14]{huppert})}
\newtheorem{thmB}{Theorem B (\cite[III.14.16]{huppert})}

\def\maxclassA{Theorem~A}
\def\maxclassB{Theorem~B}

\newcommand{\thmref}[1]{Theorem~\ref{#1}}
\newcommand{\lemref}[1]{Lemma~\ref{#1}}
\newcommand{\quotient}[2]{\raise 0.8 ex\hbox {$#1$} \big / \lower 0.8 ex
\hbox {$#2$}}

\begin{document}
\title[The Class and Coexponent of a Finite $p$-Group]
{A Bound for the Nilpotency Class of a Finite $p$-Group in terms of its 
Coexponent}

\author{Paul J. Sanders}
\address{Mathematics Institute\\ University of Warwick\\ Coventry CV4 7AL\\ UK}
\email{pjs@maths.warwick.ac.uk}
 
\author{Tom S. Wilde}
\address{Mathematics Institute\\ University of Warwick\\ Coventry CV4 7AL\\ UK}

%\begin{abstract}
%The coexponent of a finite $p$-group is introduced and we consider how the
%nilpotency class is bounded in terms of this invariant.
%\end{abstract}
 
\maketitle
%\thispagestyle{empty}
%\centerline{Mathematics Institute}
%\centerline{University of Warwick}
%\centerline{Coventry CV4 7AL, UK}

\baselineskip=18pt
\bigskip
\bigskip
\section {Introduction} 

For a prime $p$ and a finite $p$-group $P$, we define the {\bf coexponent}
of $P$ to be the least integer $f(P)$ such that $P$ possesses a cyclic
subgroup of index $p^{f(P)}$, and the purpose of this paper is to derive a
bound for the nilpotency class in terms of the coexponent. Observe that such a 
bound is possible only for odd primes since if $n > 2$, 
the group $C_{2^{n-1}} \rtimes C_2$ with the cyclic group of order 2 acting by 
inversion, is a group of order $2^n$ which has coexponent 1 and (maximal) 
nilpotency class $n-1$. We prove two theorems in this paper. 

\begin{thm}
\label{one}
Let $p$ be an odd prime and $P$ a finite \pg\ of coexponent $f(P) \geq 1$. Then 
$cl(P) \leq 2 f(P)$.
\end{thm}

From this it follows that if $P$ is a finite \pg\ with $p > 2f(P)$ then the
group can be regarded as a Lie ring by ``inverting'' the 
Baker--Campbell--Hausdorff formula (see \cite{magnus} and \cite{lazard}). This
transformation to a Lie ring setting is used in \cite{mythesis} to classify the 
finite $p$-groups of coexponent 3 for primes greater than 3. 

The bound in \thmref{one} is clearly attained by taking $P$ to be a
non-Abelian $p$-group containing a cyclic maximal subgroup (for $p > 2$), and
an examination of the proof of \thmref{one} will show that this bound is 
attained only if $f(P) = 1$ or $p = 3$. We then have

\begin{thm}
\label{two}
Let $p$ be a prime greater than 3 and $P$ a finite \pg\ of coexponent 
$f(P) \geq 2$. Then $cl(P) \leq 2 f(P) - 1$.
\end{thm}

The question of whether a 3-group $P$ exists with coexponent greater than 1 and
nilpotency class exactly $2f(P)$ is left open.
 
\section {Proofs} 

Let $p$ be an odd prime and $P$ a finite \pg\ of order $p^n$ and
coexponent $f = f(P) \geq 1$. It may be assumed that $n > 2f$ since both
Theorems are trivially true otherwise. We begin by examining the core of a
largest cyclic subgroup contained in $P$, so let $a \in P$ have order
$p^{n-f}$ and define subgroups $Q$ and $N$ of $P$ to be $\langle a \rangle$
and ${\rm Core}_P\big(\langle a \rangle\big)$ respectively.

\begin{lem}
\label{Core}
Defining integers $r$ and $s$ by
$p^r = {\rm min}\ \{\vert P : Q Q^b \vert : b~\in~P\},
\ {\rm and}\ p^{r-s} = \vert P : C_P(N) \vert$ we have
\renewcommand{\theenumi}{\roman{enumi}}
\begin{enumerate}
\item $1 \leq r \leq f$ and $\vert P : N \vert = p^{2f - r}$.
\item $s \geq 0$, and 
$[[\ldots,[N, \underbrace {P],\ldots],P}_{u\ {\rm times}}] \leq
\agemo _{(n - 2f)u} (N)$, for any integer $u \geq 1$.\\
\end{enumerate}
\end{lem}
\begin{proof}
\mbox{}
\renewcommand{\theenumi}{\roman{enumi}}
\begin{enumerate}
\item\label{lempartone}
It is easy to see that no group is the product of two 
proper conjugate cyclic subgroups and so it
follows that $1 \leq r \leq f$. To see that $\vert P : N \vert = p^{2 f - r}$
observe that for any element $b$ of $P$ we have
$$ \vert Q : Q \cap Q^b \vert = \vert QQ^b : Q \vert.$$
Now since $Q$ is cyclic there exists some element $c$
of $P$ with
${\rm Core}_P(Q) = Q \cap Q^c$ and then for any other $b \in P$ we have
$$ Q \geq Q \cap Q^b \geq Q \cap Q^c ,$$
whence $\vert QQ^c : Q \vert \geq \vert QQ^b : Q \vert $.
Therefore $\vert Q : N \vert = {\rm max}\ \{\vert QQ^b : Q \vert : b \in P \}
= p^{f - r}$ and so $\vert P : N \vert = p^{2 f - r}$.

\item\label{lemparttwo} Let $c$ be defined as in \ref{lempartone}. Then since 
$N$ is centralised 
by both $Q$ and $Q^c$ it follows that $QQ^c \leq C_P(N)$ and so $s \geq 0$.
Now let $k$ be an integer satisfying $1 \leq k \leq n - 2f +r$ (observe
by \ref{lempartone}.\ that $\vert N \vert = p^{n - 2f +r}$ and 
$n - 2f + r \geq 2$). Because $N$ is a cyclic group of prime-power
order, any
automorphism of $N / \agemo_k (N)$ lifts to an automorphism of $N$ and so
we have a composite of maps
$$P \buildrel \phi \over \longrightarrow {\rm Aut}(N)
\buildrel \gamma \over \longrightarrow 
{\rm Aut} \left ( \quotient {N}{\agemo_k (N)} \right)$$
where $P$ acts by conjugation on $N$ and $\gamma$ is onto. It follows that
$$\vert {\rm Im}(\phi) \vert = \vert P : C_P(N) \vert\ \ {\rm and}\ \  
\vert {\rm Ker}(\gamma) \vert = \frac {p^{n - 2f + r - 1} (p-1)}
{ p^{k-1} (p-1)}.$$ 
Since $p$ is odd we have that ${\rm Aut}(N)$ is cyclic, and therefore 
${\rm Im}(\phi)
\subseteq {\rm Ker}(\gamma)$ if and only if $n - 2f + r - k \geq
r - s$, i.e. if and only if $k \leq n - 2f + s$.
So taking $k = n - 2f$ we see that $[N,P] \subseteq
\agemo_k(N)$ and then the desired result follows by using induction on $u$
and the fact that for any $l \geq 0,\ \big[\agemo_l(N),P\big] \subseteq
\agemo_l\big([N,P]\big)$.
\end{enumerate}
\end{proof}

\begin{proof}[Proof of \thmref{one}]
Using the same notation as above we may assume
that $\vert P : N \vert \geq p^2$ since otherwise $P$ contains a cyclic
maximal subgroup and this is the well-known case mentioned in the introduction.
Hence if we set $k = {\rm cl}(P/N)$ then $1 \leq k \leq 2f - r - 1$
by part \ref{lempartone}.\ of \lemref{Core},
and so using part \ref{lemparttwo}.\ of \lemref{Core} we
obtain $P_{k+u+1} \subseteq \agemo_{(n- 2f)u}(N)$ for any integer $u \geq 1$.
So since $\vert N \vert = p^{n - 2f + r}$ it follows
that if $u$ is an integer greater than or equal to 1 and
$(n- 2f)u \geq n - 2f +r$ then $P_{k+u+1} = 1$.
Hence, 
\begin{eqnarray*}
{\rm cl}(P) & \leq & {\rm cl} \left ( \quotient {P}{N} \right )
+ \left \lceil \frac{n - 2f + r}{n - 2f} \right \rceil + 1 - 1 \\
& = & {\rm cl} \left ( \quotient {P}{N} \right ) + \left \lceil \frac{r}{n - 2f}
\right \rceil + 1,\\
\end{eqnarray*}
where the symbol $\lceil x \rceil$ denotes the least integer greater than or
equal to $x$ ($\in {\mathbb R}$).
Therefore, substituting for ${\rm cl}(P/N)$ we have
\begin{equation}
{\rm cl}(P) \leq 2f - r - 1 + \left \lceil \frac{r}{n - 2f} \right \rceil
+ 1 = 2f + \left \lceil r \left ( \frac{1}{n - 2f} - 1 \right ) \right \rceil
\label{Bound}
\end{equation}
which, since $n - 2f \geq 1$ by assumption,
is less than or equal to $2f$ as required.
\end{proof}

\noindent To prove \thmref{two} we determine the situations under which the 
right-hand
side of \eqref{Bound} actually attains the value $2f$, and show that unless
$f(P) = 1$ or $p = 3$ the bound on the class can be improved to $2f(P) - 1$.
As mentioned in the introduction, there exist groups with $f(P) = 1$ and 
nilpotency class 2, therefore we will assume $f(P) > 1$ so that \eqref{Bound} applies to any group we consider. We also continue to use the notation
already developed above.
Observe that there are two possible situations
under which the right-hand side of \eqref{Bound} can have the value $2f(P)$ :

\begin{enumerate}
\item\label{firstcase} $n - 2f = 1$, i.e. $\vert P \vert = p^{2f + 1}$.
\item\label{secondcase} $n - 2f > 1$ and $r = 1$.
\end{enumerate}

\noindent \lemref{caseonelem} shows that in case~\ref{firstcase}.\ the 
nilpotency class is never equal to $2f$. The proof uses the following standard
results on $p$-groups of maximal class.

\begin{thmA}
Let $G$ be a $p$-group of maximal class and order $p^n$ where 
$5\leq n\leq p+1$. Then $G/G_{n-1}$ and $G_2$ have exponent $p$.
\end{thmA}

\begin{thmB}
Let $G$ be a $p$-group of maximal class and order $p^n$ where $n>p+1$. Then
$\agemo_1(G_i) = G_{i+p-1}$ for $1\leq i\leq n-p+1$. Also, $G_1$ is a regular
$p$-group with $\Omega_1(G_1) = G_{n-p+1}$ and 
$\vert G_1/\agemo_1(G_1)\vert = p^{p-1}$.
\end{thmB}

\begin{lem}
\label{caseonelem}
Let $p$ be an odd prime and $P$ a \pg\ of coexponent
$f = f(P) > 1$ with $\vert P \vert = p^n$ where $n = 2f + 1$. Then
${\rm cl}(P) \leq 2f - 1$ {\rm(}in particular, $P$ does not have maximal class
{\rm)}.
\end{lem}

\begin{proof}
Suppose that $P$ does have maximal class $2f$ (for a
contradiction) and define $Q,N$ and $r$ as above. Since $r \geq 1$ it follows
that $\vert N \vert \geq p^2$ and because $P$ has maximal class with
$\vert P : N \vert \geq p^2$ we know that $N = P_{2f - r}$. Now,
$\vert Q \vert = p^{f+1} > p^2$ and $n = 2f + 1 \geq 5$, therefore by 
\maxclassA\ we must have $p < n-1$, in
which case we can apply \maxclassB\ to deduce that $P_{n-p+1}$ has 
exponent $p$.
Therefore $N \not \subseteq P_{n-p+1}$, i.e. $2 \leq 2f - r < n-p+1$. We can
now apply \maxclassA\ again with $i=2f-r$ to obtain 
\begin{equation}
\agemo_1(N) = P_{2f-r+p-1} \subsetneqq P_{2f-r+1} \label{Contradiction}
\end{equation}
Since $P$ has maximal class each term of the lower central series
has index $p$ in the one above (apart from $P_2$) and so we must have
$\vert N : P_{2f-r+1} \vert = p$. But since $N$ is cyclic it has a unique
subgroup of index $p$ and so $\agemo_1(N) = P_{2f-r+1}$ which contradicts
equation \eqref{Contradiction}.
\end{proof}

So we may assume that $n - 2f > 1$ and focus on case \ref{secondcase}.\ above. 
In this situation the group $P/N$ has order
$2f(P) - 1$ and contains a cyclic subgroup $Q/N$ which has index $p^{f(P)}$ and
trivial core. The next lemma shows that if $P/N$ has maximal class
and $f(P) \geq 3$
then we must have $p=3$. Thus, if we are in case~\ref{secondcase}.\ above 
with $f(P) \geq 3$ and $p > 3$ then
1 can be subtracted from the right-hand side of \eqref{Bound} when
substituting for ${\rm cl}(P/N)$ thereby bounding the class of $P$ by
$2f(P) - 1$. The proof of this lemma uses additional results on $p$-groups of
maximal class, and we indicate where they can be found in Huppert's book
\cite{huppert} as they are used.

\begin{lem}
\label{Coexbig}
Let $p$ be an odd prime and $G$ a \pg\ with
$\vert G \vert = p^{2k-1}$ where $k \geq 3$. Suppose further that $G$ contains
a cyclic subgroup $H$ of index $p^k$ which has trivial core. Then $G$ does
not have maximal class except, possibly, when $p=3$.
\end{lem}
\begin{proof}
We consider the two cases $k=3$ and $k\geq 4$ separately. 

\renewcommand{\theenumi}{\alph{enumi}}
\begin{enumerate}
\item $k=3$.\\
\noindent We suppose that $p \geq 5$ and show that $cl(G) \leq 2k - 3 = 3$, so 
let $G$ be
a group of order $p^5$ which contains a cyclic subgroup $H$ of index $p^3$
with ${\rm Core}_G(H) = 1$. Suppose (for a contradiction) that $G$ has
maximal class 4. Then since the hypotheses of \maxclassA\ are satisfied, we
know that $G / G_4$ has exponent $p$, and because $Z(G) = G_4$ we know that
$H \cap G_4 \subseteq {\rm Core}_G(H) = 1$. Therefore $HG_4 / G_4$ has order
$p^2$ and is a cyclic subgroup of $G / G_4$ which contradicts the fact that
$G / G_4$ has exponent $p$. Hence $cl(G) \leq 3$ as required. 

\item $k \geq 4$.\\
\noindent We suppose that $G$ has maximal class $2k - 2$ and show
that $p=3$. Since $n\ (=2k-1)$ is odd, $G$ is not an exceptional $p$-group of 
maximal class (by \cite[III.14.6(b)]{huppert}), i.e.\ for any $i$ with 
$2 \leq i \leq n-2$, we 
have $G_1 = C_G(G_i/G_{i+2})$ where $G_1 = C_G(G_2/G_4)$ (a proper maximal 
subgroup of $G$). Therefore by an application of 
\cite[III.14.13(b)]{huppert} it follows
that all elements of $G$ which have order greater than $p^2$ must lie in 
$G_1$. In particular, $H$ and all its conjugates are contained in 
$G_1$. So choosing $x \in G$ with $H \cap H^x = 1$ (recall that $H$ 
has trivial core)
we have that $\vert H \vert \vert H^x \vert = p^{2k-2}$. Therefore since
$G_1$ is a (proper) maximal subgroup we have $G_1 = H H^x$. Now
because the exponent of $G$ is greater than $p^2$ we must have $3 \leq p < n-1$
(by \maxclassA), and then \maxclassB\ gives us that $G_1$ is a regular
\pg. From the above
factorisation of $G_1$ we can see that $\vert G_1 : \agemo_1(G_1) \vert
\leq p^2$ and $\vert \Omega_1(G_1) \vert \geq p^2$, and so by regularity 
these two inequalities are equalities. Hence $\Omega_1(G_1) = G_{n-2}$ (since 
$\Omega_1(G_1)$ is a normal subgroup of $G$ and $G$ has maximal class). But we 
also know by \maxclassB\ that $\Omega_1(G_1) = G_{n-p+1}$, and so 
$n-p+1 = n-2$, i.e. $p=3$ as required.
\end{enumerate}
\end{proof}

We have now shown that \thmref{two} holds for all coexponents
greater than 2. Since \lemref{Coexbig} is not true for $k=2$ we deal with the
coexponent 2 case directly in the following lemma. The proof of this lemma
uses the fact that for a regular $p$-group, the terms {\it uniqueness basis} 
and {\it type invariants} make sense in direct analogy with finite Abelian 
$p$-groups. 
This result is due to Phillip Hall and the reader should consult \cite{hall} 
for the relevant details. 

\begin{lem}
Let $p$ be a prime greater than 3 and $P$ a \pg\ of
coexponent $f(P) = 2$. Then ${\rm cl}(P) \leq 2f(P) - 1$.
\end{lem}

\begin{proof}
\noindent By \thmref{one} the bound $2f(P)$ holds and so since we are assuming
$p > 2f(P)$ it follows that ${\rm cl}(P) < p$, which implies that $P$ is 
regular. Therefore, if we let $\vert P \vert = p^n$, $P$ is of $type\,(n-2,2)$ 
or $type\,(n-2,1,1)$. If $P$ is of $type\,(n-2,2)$ then $\vert P : \agemo_1(P) \vert = p^2$ 
and so $[P,P] \subseteq \agemo_1(P)$. Therefore 
$P_3 \subseteq [P,\agemo_1(P)] \subseteq \agemo_1(P_2) \subseteq \agemo_2(P)$. 
By taking
a uniqueness basis of $P$ it is straightforward to see that 
$\agemo_2(P) \subseteq Z(P)$
and therefore ${\rm cl}(P) \leq 3 = 2f(P) - 1$. If $P$ is of $type\,(n-2,1,1)$
then the $p^{\rm th}$-power of a basis element corresponding to the invariant
$n-2$ is central, and so $\vert P : Z(P) \vert \leq p^3$, from which the 
required bound follows.
\end{proof}
We have now completed the proof of \thmref{two}.

\end{document}